\documentclass[11pt,a4paper,leqno]{article}
\usepackage{amsfonts,amssymb}

\usepackage{amscd}

\vspace{0.5cm}

 4

\font\erm=cmr8

\usepackage{graphicx}

\author{M.~Dziemia\'nczuk}
\title{On multi F-nomial coefficients and Inversion formula for F-nomial coefficients}

\newtheorem{defn}{Definition}

\newtheorem{theoremn}{Theorem}
\newtheorem{observen}{Observation}

\newtheorem{corollary}{Corollary}
\newtheorem{lemma}{Lemma}

\begin{document}

\begin{center}
\noindent {\Large \textsc{On multi F-nomial coefficients and Inversion formula for F-nomial coefficients}}  \\ 

\vspace{0.5cm}

\noindent {Maciej Dziemia\'nczuk}

\vspace{0.5cm}

\noindent {\erm Student in the Institute of Computer Science, University of Gda\'nsk}

\noindent {\erm PL-80-952 Gda\'nsk, st. Wita Stwosza 57, Poland}

\noindent {\erm e-mail: Maciek.Ciupa@gmail.com}
\end{center}

\vspace{1cm}

\noindent \textbf{Summary} 

\noindent In response to \cite{akk20}, we discover the looked for inversion formula for $F$-nomial coefficients. Before supplying its proof, we generalize $F$-nomial coefficients to multi $F$-nomial coefficients and we give their combinatorial interpretation in cobweb posets language, as the number of maximal-disjoint blocks of the form $\sigma P_{k_1,k_2,...,k_s}$ of layer $\langle\Phi_1\rightarrow\Phi_n\rangle$. Then we present  inversion formula for $F$-nomial coefficients using multi F-nomial coefficients for all cobweb-admissible sequences. To this end we infer also some identities as conclusions of that inversion formula for the case of binomial, Gaussian and Fibonomial coefficients.

\vspace{0.4cm}
\noindent AMS Classification Numbers: 05A19 , 11B39, 15A09.


\vspace{0.2cm}
\noindent Presented at Gian-Carlo Polish Seminar:

\noindent \emph{http://ii.uwb.edu.pl/akk/sem/sem\_rota.htm}

\section{Preliminaries}

At first, let us recall original problem called \textbf{Inversion formula for $F$-nomial coefficients} which was brought up by A. K. Kwa\'sniewski,  in his 2001 lectures and placed then in \cite{akk20} as an Exercise 7.

\begin{quote}
	\noindent \textbf{Ex.7} \cite{akk20} 
	\emph{Discover the inversion formula i.e. the array elements 
	$\left( { n \choose k }_F \right)^{-1}$ for ${ n \choose k }_F$ 
	being the so called fibonomial coefficients, i.e.
	\begin{equation}
		{n \choose k}_F = \frac{n_F!}{k_F! (n-k)_F!},
	\end{equation}
	\noindent for $n_F = F_n$ being the $n$-th Fibonacci number ($n,k>0$).
	} (end of quote)
\end{quote}

\vspace{0.2cm}
\noindent In this note we derive inversion formula, not only for Fibonacci number, but for all F-cobweb admissible sequences \cite{akk08, akk09, akk11, md2}. Therefore we can expect general hence simpler form of new and known identities for certain sequences, such as for example for Natural and Gaussian numbers as we shall present it further on.

\section{Multi $F$-nomial coefficients}

In this section $F$-nomial coefficients ${n \choose k}_F$ and multi $F$-nomial coefficients $\break{n \choose {k_1,k_2,...,k_s}}_F$ with their respective combinatorial interpretation are considered.

\begin{defn}[\cite{akk08, akk09}]
	Let any $F$-cobweb admissible sequence, then $F$-nomial coefficients are defined as follows
	\begin{equation}
		{n \choose k}_F = \frac{n_F!}{k_F! \cdot (n-k)_F!} = \frac{n_F^{\underline{k}}}{k_F!}
	\end{equation}
	\noindent where $n_F! = n_F\cdot(n-1)_F\cdot ... \cdot 1_F$ and $n_F^{\underline{k}} = n_F\cdot(n-1)_F\cdot ...\cdot (n-k+1)_F$
\end{defn}

\noindent The combinatorial interpretation of ${n \choose k}_F$ is the following \cite{akk08, akk09}:

\vspace{0.2cm}

\begin{quote}
\noindent \emph{For $F$-cobweb tiling sequences $F$-nomial coefficient ${n \choose k}_F$ is the number of max-disjoint equipotent copies $\sigma P_m$ of the layer $\break{\langle\Phi_{k+1}\rightarrow\Phi_n\rangle}$, where $m=n-k$.
	}
\end{quote}

\vspace{0.2cm}
\noindent Now we generalize $F$-nomial to multi $F$-nomial coefficients and we give also their combinatorial interpretation in cobweb posets language.

\begin{defn}\label{def:symbol}
	Let $F\equiv\{n_F\}_{n\geq 0}$ be any natural numbers' valued sequence i.e. $n_F\in\mathbb{N}\cup\{0\}$ and $s\in\mathbb{N}$. \textbf{Multi $F$-nomial coefficient} is then identified with the symbol

\begin{equation}
	{n \choose {k_1,k_2,...,k_s}}_F = \frac{n_F!}{(k_1)_F!\cdot ... \cdot (k_s)_F!}
\end{equation}

\vspace{0.2cm}
\noindent where $k_i\in\mathbb{N}$ and $\sum_{i=1}^{s}{k_i} = n$ for $i=1,2,...,s$. In other cases is equal to zero.
\end{defn}

\begin{observen}
	Let $F$ be any $F$-cobweb admissible sequence. The value of the multi $F$-nomial coefficients is natural number or zero i.e.
	
	\begin{equation}
	{n \choose {k_1,k_2,...,k_s}}_F \in \mathbb{N} \cup \{0\}
	\end{equation}
	\noindent for any $n,k_1,k_2,...,k_s\in\mathbb{N}$.
\end{observen}

\noindent For the sake of forthcoming combinatorial interpretation we introduce the following notation.

\begin{defn}
	Let any layer $\langle\Phi_1\rightarrow\Phi_n\rangle$, $n\in\mathbb{N}$ and a composition of the number $n$ into $s$ nonzero parts designated by the vector $\langle k_1,k_2,...,k_s\rangle$, where $s\in\mathbb{N}$. Any cobweb sub-poset created by concatenate of blocks $P_{k_1},P_{k_2}, ..., P_{k_s}$ i.e.
	$$
		P_{k_1,k_2,...,k_s} = C_n\left[ F;1_F,2_F,...,(k_1)_F,1_F,...,(k_2)_F,...,1_F,...,(k_s)_F \right]
	$$
	\noindent and consequently
	$$
		\sigma P_{k_1,k_2,...,k_s} = C_n\left[ F;\sigma\langle 1_F,2_F,...,(k_1)_F,1_F,...,(k_2)_F,...,1_F,...,(k_s)_F\rangle \right]
	$$
	\noindent is called \textbf{multi-block} and denoted by $P_{k_1,k_2,...,k_s}$ and $\sigma P_{k_1,k_2,...,k_s}$.
\end{defn}

\begin{figure}[ht]
\begin{center}
	\includegraphics[width=100mm]{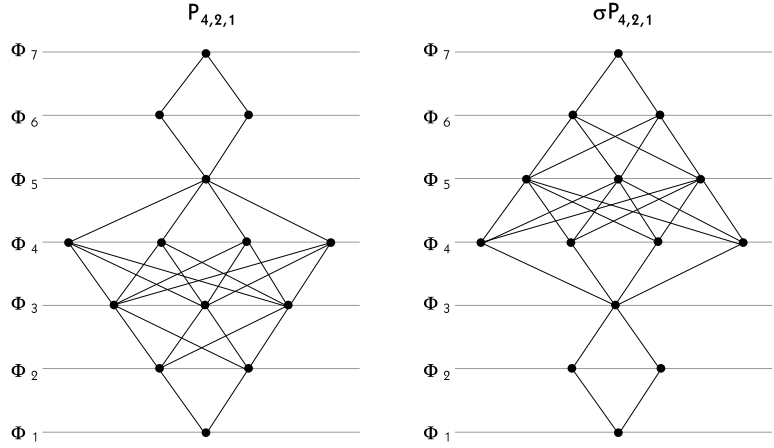}
	\caption{Picture of multi blocks $P_{4,2,1}$ and $\sigma P_{4,2,1}$. \label{fig:multiblock}}
\end{center}
\end{figure}

\begin{observen}
	For $F$-cobweb tiling sequences multi $F$-nomial coefficient $\break{n \choose {k_1,k_2,...,k_s}}_F$ is the number of max-disjoint equipotent copies $\sigma P_{k_1,k_2,...,k_s}$ \break of the layer $\langle\Phi_1\rightarrow\Phi_n\rangle$, where $n=k_1+k_2+...+k_s$.
\end{observen}

\noindent \emph{Proof.}

\noindent The number of maximal chains in a layer $\langle\Phi_1\rightarrow\Phi_n\rangle$ is equal to $n_F!$ , however the number of maximal chains in any multi block $\sigma P_{k_1,k_2,...,k_s}$ is \break $\nobreak{(k_1)_F!\cdot(k_2)_F\cdot...\cdot(k_s)_F}$. Therefore the number of blocks is equal to 

$$
	\frac{n_F!}{(k_1)_F!\cdot(k_2)_F\cdot...\cdot(k_s)_F}
$$

\vspace{0.2cm}
\noindent where $n=k_1+k_2+...+k_s$ and $n,s\in\mathbb{N}$ $\blacksquare$

\vspace{0.6cm}
\noindent Of course for $s=2$ we have

\begin{equation}
	{n \choose {k,n-k}}_F \equiv {n \choose k}_F = {n \choose {n-k}}_F
\end{equation}

\vspace{0.2cm}
\noindent Note. For any permutation $\sigma$ of the set $\{1,2,...,s\}$ the following holds

\begin{equation}
	{n \choose {k_1,k_2,...,k_s}}_F = {n \choose {k_{\sigma 1},k_{\sigma 2},...,k_{\sigma s}}}_F
\end{equation}

\vspace{0.2cm}
\noindent as is obvious from Definition \ref{def:symbol} of the multi F-nomial symbol. i.e.

$$
	\frac{n_F!}{(k_1)_F!\cdot(k_2)_F\cdot...\cdot(k_s)_F} = \frac{n_F!}{(k_{\sigma 1})_F!\cdot(k_{\sigma 2})_F\cdot...\cdot(k_{\sigma s})_F}
$$

\vspace{0.4cm}
\begin{lemma}\label{lem:1}
	Let any cobweb-tiling sequence $F$ from $\mathcal{T}_\lambda$ family \cite{md2} i.e. such that for any $m,k\in\mathbb{N}\cup\{0\}$ its terms satisfy
\begin{equation}
	n_F = (m+k)_F = \lambda_m \cdot m_F  +  \lambda_k \cdot k_F
\end{equation}

\noindent for certain coefficients $\lambda_m, \lambda_k$. Take any composition of the number $n$ into $s$ nonzero parts designated by the vector $\langle k_1,k_2,...,k_s \rangle$. Then terms of the sequence $F$ also satisfy
\begin{equation}
	n_F = \Big(\sum_{j=1}^{s}{k_j}\Big)_F = \sum_{j=1}^{s} \lambda_{k_j}\!\cdot\!(k_j)_F 
\end{equation}

\noindent for certain coefficients $\lambda_{k_j} \equiv \lambda_{k_j}(k_1,k_2,...,k_s):\mathbb{N}_0^s \rightarrow \mathbb{N}_0$, where $\nobreak{\mathbb{N}_0 \equiv \mathbb{N} \cup \{0\}}$ and $j=1,2,..,s$.
\end{lemma}

\vspace{0.2cm}
\noindent \emph{Proof.}

\noindent Let any cobweb tiling sequence $F\in\mathcal{T}_\lambda$ and take a composition of the natural number $n$ given by the vector $\langle k_1,k_2,...,k_s \rangle$. Then from definition of $\mathcal{T}_\lambda$ its terms satisfy

$$
	n_F = \left( k_1 + (n-k_1) \right)_F = \lambda_{k_1}\cdot(k_1)_F  +  \lambda_{(n-k_1)}\cdot(n-k_1)_F
$$

\vspace{0.2cm}
\noindent Then, next summand of the above can be also separates into two summands
$$
	\left(k_1 + k_2 + (n - k_1 - k_2)\right)_F = \lambda_{k_1}\cdot(k_1)_F + \lambda_{k_2}\cdot(k_2)_F + ... + \lambda_{k_s}\cdot(k_s)_F
$$

\vspace{0.2cm}
\noindent and so on up to the $k_s = (n-k_1-k_2-...-k_{s-1})$ case
$$
	(k_1+k_2+...+k_s)_F = \lambda_{k_1}\!\cdot\!(k_1)_F +  \lambda_{k_2}\!\cdot\!(k_2)_F + ... + \lambda_{k_s}\!\cdot\!(k_s)_F \ \ \blacksquare
$$

\newpage
\noindent \textbf{Examples:}

\noindent For $F$ - Natural numbers each of coefficients is constant i.e. $\lambda_{k_1},\lambda_{k_2},...,\lambda_{k_s} = const = 1$. Therefore we obtain trivial identity $(k_1+k_2+...+k_s)_F = (k_1)_F + (k_2)_F + ... + (k_s)_F$. For $F$ - Fibonacci numbers and $s=3$ case, we have
\begin{equation}
	(a+b+c)_F = \lambda_a a_F + \lambda_b b_F + \lambda_c c_F
\end{equation}

\noindent where $\lambda_a = (c+1)_F (b-1)_F$, $\lambda_b = (c+1)_F (a+1)_F$ and $\break{\lambda_a = (a_F b_F + \break{(a-1)_F (b-1)_F}}$.

\vspace{0.2cm}
\begin{theoremn}\label{th:1}
	Let any sequence $F \in \mathcal{T}_\lambda$. Then $F$ is cobweb multi tiling i.e. any layer $\langle\Phi_1\rightarrow\Phi_n\rangle$ can be tiling with the help of max-disjoint multi-blocks $\sigma P_{k_1,k_2,...,k_s}$ and the number of max-disjoint multi-blocks satisfies

\begin{equation}
	{n \choose {k_1,k_2,...,k_s}}_F = \sum_{j=1}^{s}{\lambda_{k_j}{{n-1} \choose {k_1,...,k_{j-1},k_j - 1,k_{j+1},...,k_s}}}_F
\end{equation}

\noindent where $n=k_1+k_2+...+k_s$ for any $k_1,k_2,...,k_s \in \mathbb{N}$.
\end{theoremn}

\noindent \emph{Proof.}

\noindent The main idea of this proof was already used in \cite{md5, md1}, see there for more details. Given any Cobweb poset $\Pi$ designated by sequence $F$ of cobweb tiling sequences family $\mathcal{T}_\lambda$. Let any layer $\langle\Phi_1\rightarrow\Phi_n\rangle$, where $n\in\mathbb{N}$ and any number $s\in\mathbb{N}$. We need to tiling the layer with help of max-disjoint multi blocks of the form $\sigma P_{k_1,k_2,...,k_s}$.

\vspace{0.2cm}
\noindent Consider $\Phi_n$ level, where is $n_F$ vertices, moreover from Lemma \ref{lem:1} the number of vertices in this level can be written down as the following sum

$$
	n_F = \sum_{j=1}^{s}{\lambda_{k_j}\cdot(k_j)_F}
$$

\noindent for certain coefficients $\lambda_{k_j} \equiv \lambda_{k_j}(k_1,k_2,...,k_s):\mathbb{N}^s_0 \rightarrow \mathbb{N}_0$, where $\mathbb{N}_0 \equiv \mathbb{N}\cup\{0\}$.

\begin{figure}[ht]
\begin{center}
	\includegraphics[width=80mm]{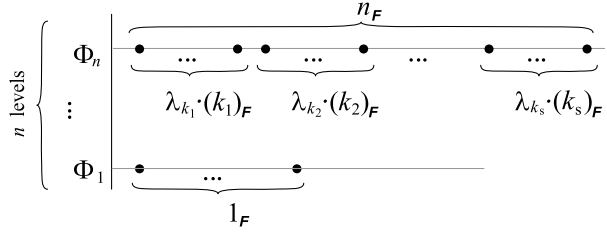}
	\caption{Idea's picture of Theorem \ref{th:1}. \label{fig:proof}}
\end{center}
\end{figure}

\vspace{0.2cm}
\noindent Therefore let us separate this $n_F$ vertices by cutting into $s$ disjoint subsets as illustrated by Fig. \ref{fig:proof} and cope at first $\lambda_{k_1}\cdot(k_1)_F$ vertices in Step 1, then $\lambda_{k_2}\cdot(k_2)_F$ ones in Step 2 and so on up to the last $\lambda_{k_s}\cdot(k_s)_F$  vertices to consider in the last $s$-th step.

\newpage
\noindent \textbf{Step 1.}

\noindent Temporarily we have $\lambda_{k_1}\cdot(k_1)_F$ fixed vertices on $\Phi_n$ level to consider. Let us cover them $\lambda_{k_1}$ times by $(k_1)$-th level of block $P_{k_1,k_2,...,k_s}$, which has exactly $(k_1)_F$ vertices. If $\lambda_{k_1} = 0$ we skip this step. What was left is the layer $\langle\Phi_1\rightarrow\Phi_{n-1}\rangle$ and we might eventually tiling it with smaller max-disjoint blocks $\sigma P_{k_1-1,k_2,...,k_s}$  in next induction step. 

\vspace{0.2cm}
\noindent \textbf{Note.} In next induction steps we use smaller blocks $\sigma P$ without level, which we have used in current step (maximal-disjoint condition). Notice also, that we could cover last level of layer $\langle\Phi_1\rightarrow\Phi_{n}\rangle$ by not the last level of block $P_{k_1,k_2,...,k_s}$ due to fact, that we can tiling by blocks obtained from levels' permutation $\sigma$ i.e. with help of maximal disjoint blocks of the form $\sigma P_{k_1,k_2,...,k_s}$.

\vspace{0.2cm}
\noindent \textbf{Step 2.}

\noindent Consider now the second situation, where we have $\lambda_{k_2}\cdot(k_2)_F$ vertices on $\Phi_n$ level being fixed. We cover them $\lambda_{k_2}$ times by $(k_1+k_2)$-th level of block $P_{k_1,k_2,...,k_s}$, which has $(k_2)_F$ vertices. Therefore we shall obtain smaller layer $\langle\Phi_1\rightarrow\Phi_{n-1}\rangle$ to tiling with blocks $\sigma P_{k_1,k_2-1,k_3,...,k_s}$.

\vspace{0.2cm}
\noindent And so on up to ...

\vspace{0.2cm}
\noindent \textbf{Step $s$.}
By analogous to previous steps, we cover the last $\lambda_{k_s}$ vertices by last $(k_1+k_2+...+k_s)=n$-th level of block $P_{k_1,k_2,...,k_s}$, obtaining smaller layer $\langle\Phi_1\rightarrow\Phi_{n-1}\rangle$ to tiling with blocks $\sigma P_{k_1,...,k_{s-1},k_s-1}$.

\vspace{0.2cm}
\noindent \textbf{Recapitulation:}

\noindent The layer $\langle\Phi_1\rightarrow\Phi_{n}\rangle$ may be tiling with blocks $\sigma P_{k_1,k_2,...,k_s}$ if  $\langle\Phi_1\rightarrow\Phi_{n-1}\rangle$ may be tiling with blocks $\sigma P_{k_1-1,k_2,...,k_s}$ and $\langle\Phi_1\rightarrow\Phi_{n-1}\rangle$ by $\sigma P_{k_1,k_2-1,k_3,...,k_s}$ again and so on up to tiling layer $\langle\Phi_1\rightarrow\Phi_{n-1}\rangle$ by $\sigma P_{k_1,...,k_{s-1},k_s-1}$. Continuing these steps by induction, we are left to prove that $\langle\Phi_1\rightarrow\Phi_{s}\rangle$  may be partitioned by blocks $\sigma P_{1,1,...,1}$ or $\langle\Phi_1\rightarrow\Phi_{1}\rangle$ by $\sigma P_{1,0,...,0}$ ones, which is obvious. $\blacksquare$

\section{Inversion formula for $F$-nomial coefficients}

Let $\Pi$ be any cobweb poset designated by $F$ - cobweb-admissible sequence. Let us consider standard reduced incidence algebra $R(\Pi)$. $F$-nomial coefficient ${n \choose k}_F$ is an element of reduced incidence algebra $R(\Pi)$ \cite{eks11,eks12}. For more, see also references therein.

\begin{theoremn}
	Let any $F$-cobweb admissible sequence and the matrix of ${n \choose k}_F$ elements. Then the inverse matrix ${n \choose k}_F^{-1}$ is given by the formula:

{\setlength\arraycolsep{2pt}
	\begin{eqnarray}
	{n \choose k}_F^{-1} & = & 
	{n \choose k}_F \bigg( \sum_{s=1}^{n-k} (-1)^s \sum_{{k_1+k_2+...+k_s=n-k} \atop {k_1,k_2,...,k_s\geq 1}} {{n-k} \choose {k_1,k_2,...,k_s}}_F \bigg) =
	\nonumber\\
	& = & \sum_{s=1}^{n-k} (-1)^s \sum_{{k_1+k_2+...+k_s=n-k} \atop {k_1,k_2,...,k_s\geq 1}} {{n} \choose {k, k_1,k_2,...,k_s}}_F
	\end{eqnarray}
}
\noindent for any $n,k \in \mathbb{N}$ such that $n\neq k$ and ${n \choose n}_F^{-1} = 1$.

\end{theoremn}

\noindent \emph{Proof.}

\noindent We need to show that ${n \choose k}_F \ast {n \choose k}_F^{-1} = \delta_{n,k}$ i.e. $\sum_{s=k}^{n}{n \choose s}_F{s \choose k}_F^{-1} = \delta_{n,k}$ for the $n\neq k$ case only, for $n=k$ it is easy to see that, it holds.

\begin{eqnarray}
	\sum_{s=k+1}^{n} {n \choose s}_F \sum_{j=1}^{s-k}(-1)^j 
	\sum_{{k_1+k_2+...+k_j=s-k} \atop {k_1,k_2,...,k_j\geq 1}} {s \choose {k,k_1,k_2,...,k_j}}_F + 
	{n \choose k}_F \cdot 1 = 0
	\nonumber\\
	\sum_{s=k+1}^{n} \sum_{j=1}^{s-k}(-1)^j
	\sum_{{k_1+k_2+...+k_j=s-k} \atop {k_1,k_2,...,k_j\geq 1}} {n \choose {k,n-s,k_1,k_2,...,k_j}}_F + 
	{n \choose k}_F = 0 
	\nonumber
\end{eqnarray}

\noindent Let us distinguish $s=n$ case from the sum:

\begin{eqnarray}
	\sum_{s=k+1}^{n-1} \sum_{j=1}^{s-k}(-1)^j
	\sum_{{k_1+k_2+...+k_j=s-k} \atop {k_1,k_2,...,k_j\geq 1}}
	{n \choose {k,n-s,k_1,k_2,...,k_j}}_F + 
	\nonumber \\ 
	 \sum_{j=1}^{n-k}(-1)^j \sum_{{k_1+k_2+...+k_j=n-k} \atop {k_1,k_2,...,k_j\geq 1}} {n \choose {k,0,k_1,k_2,...,k_j}}_F + 
	{n \choose k}_F = 0 
	\nonumber
\end{eqnarray}

\noindent Then renumerate first sum, such that $s \mapsto s - k$

\begin{eqnarray}
	\sum_{s=1}^{n-k-1} \sum_{j=1}^{s}(-1)^j
	\sum_{{k_1+k_2+...+k_j=s} \atop {k_1,k_2,...,k_j\geq 1}}
	{n \choose {k,n-s-k,k_1,k_2,...,k_j}}_F +
	\nonumber \\ 
	 \sum_{j=1}^{n-k}(-1)^j \sum_{{k_1+k_2+...+k_j=n-k} \atop {k_1,k_2,...,k_j\geq 1}} {n \choose {k,k_1,k_2,...,k_j}}_F + 
	{n \choose {k, n-k}}_F = 0 
	\nonumber
\end{eqnarray}

\noindent From the second sum we distinguish $j=1$ case and reduce it with the last summand i.e.

$$
	(-1)^1 \sum_{{k_1 = n-k} \atop {k_1\geq 1}} {n \choose {k,k_1}}_F + {n \choose {k,n-k}}_F = 0
$$

\noindent Therefore we have two summands, denote them as $A, B$ i.e. $A + B = 0$ and renumerate second sum, such that $j \mapsto j - 1$

\begin{eqnarray}
	\sum_{s=1}^{n-k-1} \sum_{j=1}^{s}(-1)^j
	\sum_{{k_1+k_2+...+k_j=s} \atop {k_1,k_2,...,k_j\geq 1}}
	{n \choose {k,n-s-k,k_1,k_2,...,k_j}}_F +
	\nonumber \\ 
	 \sum_{j=1}^{n-k-1}(-1)^{j+1} \sum_{{k_1+k_2+...+k_j+k_{j+1}=n-k} \atop {k_1,k_2,...,k_j,k_{j+1}\geq 1}} {n \choose {k,k_1,k_2,...,k_j,k_{j+1}}}_F 
	= 0
	\nonumber
\end{eqnarray}

\noindent Let us notice that, in $A$ term we have sum after Diophantine equations with non-zero terms, where $s$ takes over values from the set $\{1,2,...,n-k-1\}$ like additional variable $k_{j+1}$ in $B$ term i.e.

$$
	\left\{ 
	\begin{array}{ll}
	n = k + s + (n-k-s)
	\\
	s=1,2,...,n-k-1
	\\
	k_1+k_2+...+k_j = s
	\\
	k_1,k_2,...,k_j\geq 1
	\\
	j=1,...,s
	\end{array}
	\right. \Rightarrow
	\left\{ 
	\begin{array}{l}
	n = k + k_1 + ... + k_j + k_{j+1}
	\\
	k_1,k_2,...,k_j\geq 1
	\\
	j=1,...,s
	\end{array}
	\right.
$$

\noindent The sign of the summands in $A$ and $B$ sums changes with changing the number of variables in Diophantine equations and depends only on that, hence

\begin{eqnarray}
	\sum_{j=1}^{n-k-1} (-1)^j
	\sum_{{k_1+k_2+...+k_j+k_{j+1}=n-k} \atop {k_1,k_2,...,k_j,k_{j+1}\geq 1}}
	{n \choose {k,k_1,k_2,...,k_j,k_{j+1}}}_F +
	\nonumber \\ 
	 \sum_{j=1}^{n-k-1}(-1)^{j+1} \sum_{{k_1+k_2+...+k_j+k_{j+1}=n-k} \atop {k_1,k_2,...,k_j,k_{j+1}\geq 1}} {n \choose {k,k_1,k_2,...,k_j,k_{j+1}}}_F 
	= 0 \ \ \blacksquare
	\nonumber
\end{eqnarray}

\begin{corollary}
	Whitney numbers of first kind define inverted matrixes for Whitney numbers of second kind i.e. $\left[ W_{ij} \right]_{n\times n} \ast \left[ w_{ij} \right]_{n\times n} = \delta_{i,j}$ for certain posets \cite{lipski}. Whitney numbers of second kind of a few posets are particular case of $F$-nomial coefficients, and from Inversion formula of $F$-nomial coefficients we can infer some identities:
\end{corollary}

\begin{enumerate}
\item $\langle \mathcal{P}_n, \subseteq \rangle$ \break 
	Whitney numbers $W_{\mathcal{P}_n}(k) = {n \choose k}$ and $w_{\mathcal{P}_n}(k) = (-1)^{n-k}{n \choose k}$
	
	For $F$-\textbf{Natural numbers}
	\begin{equation}
		W_{\mathcal{P}_n}(k) = {n \choose k} = {n \choose k}_F \Rightarrow {n \choose k}_F^{-1} = (-1)^{n-k}{n \choose k}
	\end{equation}

\item $\langle L(n,q), \subseteq \rangle$ \break 
	Whitney numbers $W_{L(n,q)}(k) = {n \choose k}_q$ and $w_{L(n,q)}(k) = (-1)^{n-k}{n \choose k}_q q^{{n-k}\choose 2}$

	For $F$-\textbf{Gaussian numbers}
	\begin{equation}
		W_{L(n,q)}(k) = {n \choose k}_q\!\!\!= {n \choose k}_F\!\!\!\Rightarrow {n \choose k}_F^{-1}\!\!\!= (-1)^{n-k}{n \choose k}_q q^{{n-k}\choose 2}
	\end{equation}
\end{enumerate}

\begin{corollary}
Let $\{\Phi_n(x)\}_{n\geq 0}$ be a polynomial sequence $(\mathrm{deg}(\Phi_n(x))=n)$. Then for any $n,k \in \mathbb{N}$
\begin{equation}
	\Phi_n(x) = \sum_{k\geq 0} {n \choose k}_F^{-1} x^k
	\Leftrightarrow
	x^n = \sum_{k\geq 0} {n \choose k}_F \Phi_k(x)
\end{equation}
\noindent while $\Phi_0(x) = 1$.
\end{corollary}

\begin{enumerate}
\item \textbf{Natural numbers}, take $F$ such that $n_F = n$, then
$$
	(x - 1)^n = \sum_{k \geq 0} {n \choose k}_F (-1)^{n-k} x^k \Leftrightarrow
	x^n = \sum_{k \geq 0} {n \choose k}_F (x - 1)^k
$$
\noindent where ${n \choose k}_F \equiv {n \choose k}$
\item \textbf{Gaussian integers}, take $F$ such that $n_F = \frac{1 - q^{n-1}}{1 - q}$, then
$$
	\Phi_n(x) = \sum_{k \geq 0} {n \choose k}_F (-1)^{n-k} q^{{n-k \choose 2}} x^k \Leftrightarrow
	x^n = \sum_{k \geq 0} {n \choose k}_F \Phi_k(x)
$$
\noindent where ${n \choose k}_F \equiv {n \choose k}_q$ and $\Phi_n(x) = \prod_{s=0}^{n-1} (x - q^s)$.
\item \textbf{Fibonacci numbers}, i.e. $n_F = (n-1)_F + (n-2)_F$, $1_F=2_F=1$ then with the help of inversion formula for F-nomial coefficients, let us show a few first polynomials $\Phi_n(x)$ 

\noindent $\Phi_0(x) = 1$

\noindent $\Phi_1(x) = x - 1$

\noindent $\Phi_2(x) = x^2 - x$

\noindent $\Phi_3(x) = x^3 - 2x^2 + 1$

\noindent $\Phi_4(x) = x^4 - 3x^3 + 3x - 1$

\noindent $\Phi_5(x) = x^5 - 5x^4 + 15x^2 - 5x - 6$

\noindent $\Phi_6(x) = x^6 - 8x^5 + 60x^3 - 40x^2 - 48x + 35$

\noindent $\Phi_7(x) = x^7 - 13x^6 + 260x^4 - 260x^3 - 624x^2 + 455x + 181$.

\noindent $\Phi_8(x) = x^8 - 21x^7 + 1092x^5 - 1820x^4 - 6552x^3 + 9555x^2 + 8301x - 6056$.

\noindent however simpler formula for $\Phi_n(x)$ where $F$-Fibonacci sequence as for instance in the case of Natural and Gaussian numbers is still unknown.
\end{enumerate}

\noindent \textbf{To be the next}

\noindent As we can see above, for certain sequences $F$, inversion formula reduces to simpler form. Therefore for other sequences, like Fibonacci numbers with a lot of properties, we can expect also simpler ones.

\vspace{0.4cm}
\noindent \textbf{Acknowledgements}

\vspace{0.4cm}
I would like to thank Professor A. Krzysztof Kwa\'sniewski - who initiated my interest in his cobweb poset concept - for his very helpful comments, improvements and corrections of this note. Assistance of E. Krot-Sieniawska with respect to incidence algebras of cobweb posets is highly appreciated too.


\end{document}